\documentclass{article}

\usepackage{amssymb,amscd,amsmath,amsxtra}
\usepackage[dvips]{graphicx}
\usepackage[arrow,matrix,curve]{xy}

\newtheorem{Thm}{Theorem}
\newtheorem{Prop}[Thm]{Proposition}
\newtheorem{Lem}[Thm]{Lemma}

\newtheorem{Cor}[Thm]{Corollary}
\newcommand{\eop}{\textrm{\hspace{1em}$\square$}}

\DeclareMathOperator{\Rep}{Rep}

\newcommand{\bZ}{\mathbb Z}
\newcommand{\bF}{\mathbb F}
\renewcommand{\a}{\mathbf a}
\newcommand{\sI}{\mathcal I}
\newcommand{\bI}{\mathbf I}
\newcommand{\bi}{\mathbf i}
\newcommand{\bj}{\mathbf j}

\begin{document}

\begin{center}
\large Quiver representations respecting a quiver automorphism:\\
a generalisation of a theorem of Kac
\end{center}

\begin{center}
Andrew Hubery\\
University of Leeds, England\\
ahubery@maths.leeds.ac.uk
\end{center}

\begin{abstract}
A theorem of Kac on quiver representations states that the dimension vectors of indecomposable representations are precisely the positive roots of the associated symmetric Kac-Moody Lie algebra. Here we generalise this result to representations respecting an admissible quiver automorphism, and  obtain the positive roots of an associated symmetrisable Kac-Moody Lie algebra.

We also discuss the relationship with species of valued quivers over finite fields. It is known that the number of isomorphism classes of indecomposable representations of a given dimension vector for a species is a polynomial in the size of the base field. We show that these polynomials are non-zero if and only if the dimension vector is a positive root of the corresponding symmetrisable Kac-Moody Lie algebra.
\end{abstract}

MSC (2000): 16G20 (17B67 16S35).

\section{Introduction}

A remarkable theorem by Kac \cite{Kac1} states that the dimension vectors of indecomposable representations of a quiver $Q$ without vertex loops over an algebraically closed field are precisely the positive roots of an associated Kac-Moody Lie algebra $\mathfrak{g}$. The Kac-Moody Lie algebras that appear in this way are the symmetric ones --- i.e. those with a symmetric generalised Cartan matrix (GCM). However, there is also a well-developed theory for the symmetrisable Kac-Moody algebras, whose GCM can be symmetrised by left multiplication by an invertible diagonal matrix. It would clearly be of interest to generalise Kac's Theorem in order to obtain the positive root systems of these symmetrisable Kac-Moody algebras.

We know that from any symmetrisable GCM, with valued graph $\Gamma$, we can construct a graph with an admissible automorphism $(Q,\a)$ \cite{Lu}. This pair is in no way unique, but there is a close connection between properties of $Q$ compatible with $\a$ and the corresponding properties of $\Gamma$. For example, there is an obvious bijection
$$f:(\bZ\sI)^{\langle\a\rangle}\to\bZ\bI$$
between the fixed points of the root lattice for $Q$ and the root lattice for $\Gamma$, which in turn gives rise to an isomorphism between the Weyl group $W(\Gamma)$ and the subgroup of $W(Q)$ of elements commuting with $\a$. (Here, $\sI$ and $\bI$ are indexing sets for the vertices of $Q$ and $\Gamma$ respectively.)

The automorphism $\a$ naturally acts on the category of representations of $Q$ as an autoequivalence. It is therefore interesting to see how it acts on the isomorphism classes of indecomposable representations. This action is perhaps best understood by considering the following representations. We take one representative from each isomorphism class in the orbit and form their direct sum. This representation is then isomorphic to its twist by $\a$, and we refer to such representations as isomorphically invariant indecomposables (or ii-indecomposables). The dimension vector $\alpha$ of an ii-indecomposable is clearly fixed by $\a$ and so we can consider the element $f(\alpha)$ in the root lattice for $\Gamma$.

We can now state the main result of this paper.
\begin{Thm}\label{MainThm}
Let $(Q,\a)$ be a quiver with admissible automorphism and let $k$ be an algebraically closed field of characteristic not dividing the order of $\a$. Then
\begin{enumerate}
\item the images under $f$ of the dimension vectors of the ii-indecomposables give the positive root system for the symmetrisable Kac-Moody algebra  $\mathfrak{g}(\Gamma)$;
\item if $f(\alpha)$ is real, there is a unique such ii-indecomposable of dimension vector  $\alpha$ up to isomorphism. The number of direct summands of this ii-indecomposable equals the length of the root $f(\alpha)$. (Equivalently, this is the size of the corresponding orbit of isomorphism classes of indecomposables.)
\end{enumerate}
\end{Thm}
Here, the length of a root $\beta$ for $\Gamma$ indicates the value $\tfrac{1}{2}(\beta,\beta)_\Gamma$ with respect to the bilinear form (see \S2).

In particular, the dimension vectors of the ii-indecomposables are independent of the (compatible) orientation of $Q$ and of the characteristic of the field $k$ (assuming it does not divide the order of $\a$).

This extends a result of Tanisaki \cite{Tan} which states that the pair $(Q,\a)$ has only finitely many isomorphism classes of ii-indecomposables if and only if $Q$ is Dynkin.

In fact, we can consider the following more general situation. We know that the category of representations of a quiver $Q$ over a field $k$ is equivalent to the category of finite dimensional modules for the path algebra $kQ$ \cite{ARS}. This path algebra can be thought of as a tensor algebra
$T(\Lambda_0,\Lambda_1)$, where $\Lambda_0=\prod_{i\in\sI}k\varepsilon_i$ is a basic semi-simple $k$-algebra and $\Lambda_1=\coprod_{\rho:i\to j}k\rho$ is a finite dimensional $\Lambda_0$-bimodule. To say that $Q$ has no edge loops is then equivalent to specifying that $\varepsilon_i\Lambda_1\varepsilon_i=0$ for all $i$.

Now suppose that $\a$ is any algebra automorphism of $T(\Lambda_0,\Lambda_1)$ given by an algebra automorphism of $\Lambda_0$ and a bimodule automorphism of $\Lambda_1$ (so $\a$ respects the natural grading of $T$). This induces an action of $\a$ on the vertex set $\sI$, and we say that this action is admissible if $\varepsilon_j\Lambda_1\varepsilon_i=0$ whenever $i$ and $j$ lie in the same $\a$-orbit.

We can again associate to this pair $(Q,\a)$ a symmetrisable GCM, and hence a valued graph $\Gamma$, and all the results of this paper extend without alteration to this setting. (In fact, we are forced to consider such automorphisms when we consider the skew group algebra $kQ\#\langle\a\rangle$.)

The heart of the proof of the main theorem lies in understanding the relationship between ii-indecomposables for $(Q,\a)$ and indecomposables for the skew group algebra $kQ\#\langle\a\rangle$, the main reference being the paper by Reiten and Riedtmann \cite{RR}.

In the penultimate section, we describe some particular examples of quivers with an automorphism. We note that if $Q$ has no oriented cycles, then the action of $\a$ on the category of representations commutes with that of the Auslander-Reiten translate $\tau$. In particular, if $Q$ is a tame quiver, then the actions of $\a$ and $\tau$ on the regular modules are determined by their actions on the regular simples, and hence the ii-indecomposables can be classified. 

Finally we explain the correspondence between ii-representations and modules for the $k$-species associated to $\Gamma$. In particular, let $k$ be a finite field and $t$ the highest common factor of the so-called symmetrisers of $\Gamma$ (c.f. \cite{Hua2}). Write $\Lambda$ for the tensor algebra associated to the $k$-species of $\Gamma$. Then for $K/k$ an extension of degree $t$, the $K$-algebra $K\otimes\Lambda$ is isomorphic to the path algebra of a quiver $KQ$ on which the Galois group $\mathrm{Gal}(K,k)$ acts as $K$-algebra automorphisms. It turns out that this is actually generated by an admissible quiver automorphism $\a$ of $Q$, and the valued graph associated to the pair $(Q,\a)$ is the underlying graph of $\Gamma$. Moreover, we have the standard Galois group action (as $k$-algebra automorphisms) coming from the identification $KQ\cong K\otimes_kkQ$, and the induced modules $K\otimes X$ for $X$ a $\Lambda$-module are precisely those for which these two actions coincide.

We know from \cite{Hua2} that the number of isomorphism classes of indecomposable $\Lambda$-modules of dimesion vector $\alpha$ is a polynomial in the size of the base field $k$. Here we show that these polynomials are non-zero precisely when $\alpha$ is a positive root of $\Gamma$, thus answering a question of Hua \cite{Hua2}.

This paper is based on work done during the author's doctorate at the University of Leeds, supported by the Engineering and Physical Sciences Research Council. The author would like to thank his supervisor Prof.~W.~Crawley-Boevey for his continuing support and guidance. Part of this paper was written whilst the author was visiting BUGH Wuppertal, supported by the ESF/PESC Programme on Non-Commutative Geometry, and he would like to thank Dr.~M.~Reineke, Prof.~Dr.~K.~Bongartz and the other members of the algebra group for their hospitality.

\section{Quivers With An Admissible Automorphism}

In this paper we shall only consider finite quivers without vertex loops. Given such a quiver $Q$, we shall denote its vertex set by $\sI$. The  corresponding symmetric GCM is the matrix $A$ indexed by $\sI$ with entries
$$a_{ij}:=\begin{cases} 2 &\textrm{if } i=j;\\
-\#\{\textrm{edges between vertices $i$ and $j$}\} &\textrm{if } i\neq j.\end{cases}$$
We write $\mathfrak{g}(Q)$ for the associated symmetric Kac-Moody algebra, with root system $\Delta(Q)$ \cite{Kac4}.

In \cite{Kac1} (see also \cite{Kac3, KR}) Kac proved the following theorem.
\begin{Thm}[Kac]
Let $Q$ be a quiver and $k$ an algebraically closed field. Then
\begin{enumerate}
\item the dimension vectors of the indecomposable representations of $Q$ over $k$ are precisely the positive roots of the Kac-Moody algebra
$\mathfrak{g}(Q)$;
\item there exists a unique isomorphism class of indecomposables of  dimension vector $\alpha$ if and only if $\alpha$ is a positive real root of $\mathfrak{g}(Q)$.
\end{enumerate}
\end{Thm}

As a consequence, the dimension vectors of the indecomposable representations are independent of both the orientation of $Q$ and the characteristic of the field $k$.

An admissible automorphism $\a$ of a quiver $Q$ is a quiver automorphism such that no arrow connects two vertices in the same orbit. This notion was first introduced by Lusztig in \cite{Lu}, where given such a pair $(Q,\a)$, he describes how to construct a symmetric matrix $B$ indexed by the vertex orbits $\bI$. Namely
$$b_{\bi\bj}:=\begin{cases} 2\#\{\textrm{vertices in \textbf{i}-th orbit}\} &\textrm{if }\bi=\bj;\\ -\#\{\textrm{edges between \textbf{i}-th and \textbf{j}-th orbits}\} &\textrm{if }\bi\neq\bj.\end{cases}$$
Let
$$d_\bi:=b_{\bi\bi}/2=\#\{\textrm{vertices in \textbf{i}-th orbit}\}$$
and set $D=\mathrm{diag}(d_\bi)$. Then $C=D^{-1}B$ is a symmetrisable GCM and we write $\Gamma$ for the corresponding valued graph. That is, $\Gamma$ has vertex set $\bI$ and we draw an edge $\bi-\bj$ equipped with the ordered pair $(|c_{\bj\bi}|,|c_{\bi\bj}|)$ whenever $c_{\bi\bj}\neq0$.

N.B. Since $\a$ is admissible, $\Gamma$ has no vertex loops.

For example, consider the quiver $A_3$ with automorphism
$\quad\vcenter{\xymatrix@R=5pt{\cdot\ar@/_/@{<.>}[dd]\ar[dr]\\&\cdot\\\cdot\ar[ur]}}$
The corresponding symmetrisable GCM is then
$(\begin{smallmatrix}2&-1\\-2&2\end{smallmatrix})=(\begin{smallmatrix}2&0\\0&1\end{smallmatrix})^{-1}(\begin{smallmatrix}4&-2\\-2&2\end{smallmatrix})$
and the valued graph $\Gamma$ is $\xymatrix@1{\cdot\ar@{-}[r]^{(2,1)}&\cdot}$

The roots of $\Gamma$ can be defined combinatorially as follows. For a vertex $\bi$, we define the reflection $r_\bi$ on the root lattice $\bZ\bI$ by
$$r_\bi:\alpha\mapsto \alpha-\frac{1}{d_\bi}(\alpha,e_\bi)_\Gamma e_\bi,$$
where $(-,-)_\Gamma$ is the symmetric bilinear form determined by the matrix $B$ and the $e_\bi$ form the standard basis of $\bZ\bI$. The group generated by these reflections is the Weyl group $W(\Gamma)$. The real roots are given by the images under $W(\Gamma)$ of the simple roots $e_\bi$ and the imaginary roots are given by $\pm$ the images under $W(\Gamma)$ of the fundamental roots
$$F_\Gamma:=\{\alpha>0\mid(\alpha,e_\bi)_\Gamma \leq 0 \textrm{ for all }\bi\textrm{ and supp}(\alpha)\textrm{ connected}\}.$$
Similarly, we have the Weyl group and the roots associated to the quiver $Q$, this time using the symmetric matrix $A$ to define $(-,-)_Q$. That is, the reflection $r_i$ on $\bZ\sI$ is given by $\alpha\mapsto\alpha-(\alpha,e_i)_Qe_i$.

The automorphism $\a$ acts naturally on the root lattice $\bZ\sI$ for $Q$, and the bilinear form $(-,-)_Q$ is $\a$-invariant. We have a canonical bijection
$$f:(\bZ\sI)^{\langle\a\rangle}\to\bZ\bI$$
from the fixed points in the root lattice for $Q$ to the root lattice for $\Gamma$. This is given by $f(\alpha)_\bi:=\alpha_i$ for any vertex $i$ in the $\bi$-th orbit.

The admissibility of $\a$ implies that the reflections $r_i$ and $r_j$ commute whenever $i$ and $j$ lie in the same $\a$-orbit. Therefore the element
$$s_\bi:=\prod_{i\in\bi}r_i\in W(Q)$$
is well-defined. Also, since $\a\cdot r_i=r_{\a(i)}\cdot\a$, we see that $s_\bi\in C_\a(W(Q))$, the set of elements in the Weyl group commuting with the action of $\a$.
\begin{Lem}\label{Lem} For $\alpha,\beta\in(\bZ\sI)^{\langle\a\rangle}$ we have
\begin{enumerate}
\item $(\alpha,\beta)_Q=(f(\alpha),f(\beta))_\Gamma$;
\item $f(s_\bi(\alpha))=r_\bi(f(\alpha))\in\bZ\bI$;
\item the map $r_\bi\mapsto s_\bi$ induces an isomorphism $W(\Gamma)\xrightarrow{\sim}C_\a(W(Q))$.
\end{enumerate}
\end{Lem}

\textit{Proof.} Let $i_v\in\sI$ enumerate the vertices in the orbit $\bi\in\bI$. Parts 1. and 2. now follow from the formula
$$b_{\bi\bj}=\sum_{v,w}a_{i_vj_w}=d_\bi\sum_w a_{i_vj_w}\quad\textrm{for any }v.$$

We denote the length of an element $w\in W(Q)$ by $\ell(w)$. Then
$$\ell(wr_i)<\ell(w) \quad\textrm{if and only if}\quad w(e_i)<0$$
and
$$\ell(w)=\#\{\alpha\in\Delta(Q)_+\mid w(\alpha)<0\}.$$
Since $\a$ is admissible and preserves the partial order $<$ on $\bZ\sI$, induction on length shows that $C_\a(W(Q))$ is generated by the $s_\bi$.
That $r_\bi$ and  $s_\bi$ satisfy the same relations follows from 2. \eop

\begin{Prop}\label{roots}
Let $\alpha\in\Delta(Q)$ and let $r\geq1$ be minimal such that $\a^r(\alpha)=\alpha$. Set
$$\sigma(\alpha):=\alpha+\a(\alpha)+\cdots+\a^{r-1}(\alpha)\in(\bZ\sI)^{\langle\a\rangle}.$$
Then $\alpha\mapsto f(\sigma(\alpha))$ induces a surjection $\Delta(Q)\twoheadrightarrow\Delta(\Gamma)$. Moreover, if $f(\sigma(\alpha))$ is real, then $\alpha$ must also be real and unique up to $\a$-orbit.
\end{Prop}

\textit{Proof.} Set $\beta:=f(\sigma(\alpha))$ and consider $w'(\beta)$ for some $w'\in W(\Gamma)$. Let $w\in C_\a(Q)$ correspond to $w'$. Then $w(\sigma(\alpha))=\sigma(w(\alpha))$ and $w(\alpha)\in\Delta(Q)$.

Since $\a$ preserves the partial order $<$ on $\bZ\sI$, $w'(\beta)$ is either positive or negative. Also, $f(\sigma(w(\alpha)))$ always has connected support. Therefore, if $w'$ is chosen so that $w'(\beta)$ has minimal height, then either $w'(\beta)$ lies in the fundamental region or else is a multiple of a real root, say $w'(\beta)=me_\bi$. In the latter case we must have that $w(\alpha)=me_i$ for some vertex $i$ in the $\bi$-th orbit, so $w(\alpha)=e_i$ and $m=1$. This proves that every root of $Q$ gives rise to a root of $\Gamma$.

Conversely, any simple root $e_\bi$ clearly lies in the image of this map, so let $\beta$ lie in the fundamental region for $\Gamma$. Then $\gamma=f^{-1}(\beta)\in(\bZ\sI)^{\langle\a\rangle}$ satisfies
$$0\geq (\beta,e_\bi)_\Gamma=(\gamma,\sigma(e_i))_Q=\sum_r(\gamma,\a^r(e_i))_Q=\frac{1}{d_\bi}(\gamma,e_i)_Q$$
for all $i$, where $i\in\bi$. Thus any connected component $\alpha$ of $\gamma$ lies in the fundamental region for $Q$ and $\sigma(\alpha)=\gamma$.

The result now follows from the correspondence between the Weyl groups. \eop

\section{Isomorphically Invariant Representations}

Let $X$ be a representation of $Q$, given by vector spaces $X_i$ for $i\in\sI$ and linear maps $X_\rho:X_i\to X_j$ for each arrow $\rho:i\to j$. The dimension vector $\underline\dim\,X\in\bZ\sI$ is defined to be $\sum_i(\dim X_i)e_i$.

We define a new representation ${}^\a X$ by $({}^\a X)_i:=X_{\a^{-1}(i)}$ and $({}^\a X)_\rho:=X_{\a^{-1}(\rho)}$. Similarly, given a homomorphism $\phi:X\to Y$, we have ${}^\a\phi:{}^\a X\to {}^\a Y$ given by $({}^\a\phi)_i:=\phi_{\a^{-1}(i)}:({}^\a X)_i\to({}^\a Y)_i$. This determines a functor $F(\a)$ on the category of representations. Clearly $F(\a^r)=F(\a)^r$ and $F(\a)$ is additive. In fact, $F(\a)$ is an autoequivalence of the category, so in particular $X$ is indecomposable if and only if ${}^\a X$ is.

We say that $X$ is isomorphically invariant (an ii-representation) if ${}^\a X\cong X$. Note that $\underline\dim\, {}^\a X=\a(\underline\dim\, X)$, so any ii-representation has dimension vector fixed by $\a$. We say that $X$ is an ii-indecomposable if it is not isomorphic to the proper direct sum of two such ii-representations. By the Krull-Remak-Schmidt Theorem for $Q$-representations, the ii-indecomposables $X$ are precisely the representations of the form
$$X\cong Y\oplus {}^\a Y\oplus\cdots\oplus {}^{\a^{r-1}}Y,$$
where $Y$ is an indecomposable representation and $r\geq 1$ is minimal such that ${}^{\a^r}Y\cong Y$.

\begin{Lem}
The Krull-Remak-Schmidt Theorem holds for ii-representations.
\end{Lem}

\textit{Proof.} Suppose that $X$ is an ii-representation. We can write $X$ as a direct sum of indecomposable representations and, since ${}^\a X\cong X$, $F(\a)$ must act (up to isomorphism) as a permutation of these indecomposable summands. That is, we can write
$$X\cong Y_1\oplus\cdots\oplus Y_m$$
with the $Y_j$ ii-indecomposables. The uniqueness properties now follow from the Krull-Remak-Schmidt Theorem for  $Q$-representations. \eop

\begin{Prop}
Let $(Q,\a)$ be a quiver with an admissible automorphism, $\Gamma$ the associated valued graph and $k$ an arbitrary algebraically closed field. Then there is an ii-indecomposable of dimension vector $\alpha$ only if $f(\alpha)$ is a root of $\Gamma$. Moreover, every positive real root of $\Gamma$ occurs, and the corresponding ii-indecomposable is unique up to isomorphism with $\frac{1}{2}(\alpha,\alpha)_Q$ indecomposable summands.
\end{Prop}

\textit{Proof.} Let $Y$ be an indecomposable representation of $Q$ and consider the ii-indecomposable
$$X\cong Y\oplus {}^\a Y\oplus\cdots\oplus {}^{\a^{r-1}}Y$$
Writing $\beta=\underline\dim\, Y\in\Delta(Q)_+$, then $\underline\dim\, X=\beta+\a(\beta)+\cdots+\a^{r-1}(\beta)=m\sigma(\beta)$ for some $m$. If $\beta$ is real, then by Kac's Theorem $Y\cong {}^{\a^r}Y$ if and only if $\a^r(\beta)=\beta$ and so $m=1$. Therefore $f(\underline\dim\, X)=f(\sigma(\beta))\in\Delta(\Gamma)_+$. On the other hand, if $\beta$ is imaginary, then so is $f(\sigma(\beta))$, and hence $f(\underline\dim\, X)=mf(\sigma(\beta))\in\Delta(\Gamma)_+$.

This shows that the dimension vectors of ii-indecomposables give rise to positive roots of $\Gamma$. Also, it follows from Proposition \ref{roots} that we must get every real root of $\Gamma$, and that the corresponding ii-indecomposable is unique up to isomorphism with the stated number of indecomposable summands. \eop

Of special interest in the representation theory of quivers are the reflection functors $R_i^+$ and $R_i^-$, defined when $i\in\sI$ is a sink or a source respectively. Clearly if $i$ is a sink (respectively a source) then the same is true for all vertices in the orbit of $i$. Also, since $\a$ is admissible, the functor
$$S_\bi^\pm:=\prod_{i\in\bi} R_i^\pm$$
is well-defined.

Denote by $k_i$ the $i$-th simple representation of $Q$. If $d=d_\bi$ is the size of the orbit of $i$, then we have an ii-indecomposable
$$k_\bi:=k_i\oplus k_{\a(i)}\oplus\cdots k_{\a^{d-1}(i)}.$$

\begin{Prop}
Let $i$ be a sink (or a source, interchanging $+$ and $-$). Then for any ii-representation $X$ there is a canonical monomorphism
$$\phi_X:S_\bi^-S_\bi^+(X)\to X$$
whose image has a complement a direct sum of copies of the ii-indecomposable $k_\bi$. In fact,
\begin{enumerate}
\item $S_\bi^+(k_\bi)=0$;
\item if $X\neq k_\bi$ is an ii-indecomposable, then $\phi_X$ is an isomorphism and hence
$$\mathrm{End}(S_\bi^+(X))\cong\mathrm{End}(X) \qquad\textrm{and}\qquad \underline\dim\, S_\bi^+(X)=s_\bi(\underline\dim\, X).$$
\end{enumerate}
\end{Prop}

\textit{Proof.} This follows immediately from the standard properties of reflection functors. \eop

\section{Skew Group Algebras}

We recall \cite{ARS} that the category of representations of a quiver $Q$ over a field $k$ is equivalent to the category of finite dimensional modules for the path algebra $kQ$. If $\a$ is a quiver automorphism of $Q$, then it naturally induces an algebra automorphism on $kQ$, also denoted $\a$. In
such a situation, we can form the skew group algebra $kQ\#\langle\a\rangle$. This has as basis the elements $\lambda\a^r$, with $\lambda$ a path in $kQ$, and the multiplication is given by
$$\lambda\a^r\cdot\mu\a^s:=\lambda\a^r(\mu)\a^{r+s}.$$

The action of $\a$ on $\bmod\, kQ$ is given by $\mathcal X\mapsto{}^\a\mathcal X$, where ${}^\a\mathcal X$ has the same underlying vector space as $\mathcal X$, but with the new action
$$\lambda\cdot x:=\a^{-1}(\lambda)x \quad\textrm{for all }\lambda\in kQ.$$

Now suppose that $k$ is an algebraically closed field with characteristic not dividing the order $n$ of $\a$. Then, given any ii-representation $X$, there exists an isomorphism $\theta:{}^\a X\to X$ and, as observed in \cite{Gab}, this implies the existence of an isomorphism $\phi:{}^\a  X\to X$
such that
$$\phi\,{}^\a\phi\,\cdots\,{}^{\a^{n-1}}\phi=1.$$
In terms of modules for the path algebra $kQ$, this says that we can find a vector space isomorphism $\Phi:\mathcal{X}\to\mathcal{X}$ such that $\Phi^n=1$ and $\Phi(\lambda x)=\a(\lambda)\Phi(x)$. That is, the pair $(\mathcal{X},\Phi)$ determines a module for the skew group algebra $kQ\#\langle\a\rangle$. We note that there may exist $\Phi$ and $\Phi'$ such that $(\mathcal X,\Phi)$ and $(\mathcal X,\Phi')$ are non-isomorphic.

We now consider a slightly more general situation, following the paper by Reiten and Riedtmann \cite{RR}. Let $k$ be an algebraically closed field and $G$ a finite cyclic group whose order is invertible in $k$. Let $\Lambda$ be either an artin $k$-algebra or the path algebra of a quiver and suppose that $G$ acts as either algebra automorphisms or graded algebra automorphisms of $\Lambda$. We write $\Lambda G$ for the skew group algebra. The following are proved in \cite{RR} when $\Lambda$ is artinian, but the proofs carry over to the path algebra case also.

Consider the following pair of functors. These are both left and right adjoints of each other (\cite{RR}, Theorem 1.1).
\begin{alignat*}{2}
F&:=\Lambda G\otimes_\Lambda- &&:\bmod\Lambda \to \bmod\Lambda G\\
H&:=\textrm{restriction} &&:\bmod\Lambda G \to \bmod\Lambda.
\end{alignat*}
We have (\cite{RR}, Proposition 1.8)
\begin{Prop}\label{Prop1}
Let $X,Y$ be indecomposable $\Lambda$-modules. Then
\begin{enumerate}
\item $HF(X)\cong\bigoplus_{g\in G}{}^g X$;
\item $F(X)\cong F(Y)$ if and only if $Y\cong {}^g X$ for some $g\in G$;
\item $F(X)$ has exactly $m$ indecomposable summands, where $m$ is the order of $\{g\in G\mid {}^gX\cong X\}$.
\end{enumerate}
\end{Prop}

It is also shown that the dual group $\hat G$ acts on $\Lambda G$ via $\chi(\lambda g)=\chi(g)\lambda g$ and then (\cite{RR}, Proposition 5.1)
\begin{Prop}
The map $\phi:(\Lambda G)\hat G\to\mathrm{End}_\Lambda(\Lambda G)$ given by
$$\phi(\lambda g\chi)(\mu h):=\chi(h)\lambda g\mu h$$
is an algebra isomorphism.

In particular, since $\Lambda G$ is a finitely generated projective generator for $\bmod\,\Lambda$, we have that $(\Lambda G)\hat G$ is Morita equivalent to $\Lambda$.
\end{Prop}

We can also consider the functors
\begin{alignat*}{2}
F'&:=(\Lambda G)\hat G\otimes_{\Lambda G}- &&:\bmod\Lambda G \to \bmod(\Lambda G)\hat G\\
H'&:=\textrm{restriction} &&:\bmod(\Lambda G)\hat G \to \bmod\Lambda G,
\end{alignat*}
and we note that the Morita equivalence $\bmod(\Lambda G)\hat G \to \bmod\Lambda$ is given by
$$M:=H\circ\Lambda G\otimes_{(\Lambda G)\hat G}-.$$
For, this comes from viewing $\Lambda G$ as a finitely generated projective generator for $\bmod\,\Lambda$.

\begin{Cor}\label{Cor1}
There are natural isomorphisms
$$H\cong MF' \qquad\textrm{and}\qquad H'\cong FM.$$
\end{Cor}

\textit{Proof.} It is clear that $MF'$ is naturally isomorphic to $H$. The second isomorphism follows by taking adjoints. \eop

It follows that Proposition \ref{Prop1} holds with $F$ and $H$ interchanged. Namely
\begin{Prop}
Let $X,Y$ be indecomposable $\Lambda G$-modules. Then
\begin{enumerate}
\item $FH(X)\cong\bigoplus_{\chi\in\hat G}{}^\chi X$;
\item $H(X)\cong H(Y)$ if and only if $Y\cong {}^\chi X$ for some $\chi\in\hat G$;
\item $H(X)$ has exactly $m$ indecomposable summands, where $m$ is the order of $\{\chi\in\hat G\mid {}^\chi X\cong X\}$.
\end{enumerate}
\end{Prop}

\begin{Lem}\label{Lem2}
Let $e$ be an idempotent of $\Lambda$ such that $e\Lambda e$ is Morita equivalent to $\Lambda$ and suppose that $e$ is fixed by $G$. Then we have natural isomorphisms between
\begin{align*}
{}&\xymatrix@C=50pt{\bmod\Lambda\ar[r]^{\Lambda G\otimes_\Lambda-}&\bmod\Lambda G\ar[r]^{e\Lambda G\otimes_{\Lambda G}-}&\bmod\, e\Lambda Ge}\\
and\qquad&\xymatrix@C=50pt{\bmod\Lambda\ar[r]^{e\Lambda\otimes_\Lambda-}&\bmod\, e\Lambda e\ar[r]^{e\Lambda Ge\otimes_{e\Lambda e}-}&\bmod\, e\Lambda Ge}\qquad
\end{align*}
\end{Lem}

We can now relate this back to quivers with an admissible automorphism.

Let $(Q,\a)$ be a quiver with an admissible automorphism of order $n$ and let $k$ be an algebraically closed field of characteristic not dividing $n$. Then the skew group algebra $kQ\#\langle\a\rangle$ is Morita equivalent to the path algebra of another quiver $\widetilde Q$ and if $\tilde\a$ is a generator for the dual group action on $kQ\#\langle\a\rangle$, then we have an induced action on $k\widetilde Q$. N.B. This may not come from a quiver automorphism, but is admissible in the sense that no edge joins two vertices in the same orbit (see the introduction).

Therefore we can consider the induction and restriction functors $F$ and $H$ as going between $\Rep(Q)$ and $\Rep(\widetilde Q)$.

On the other hand, we can do the same thing starting with $(\widetilde Q,\tilde\a)$. In this way we recover our original pair $(Q,\a)$. Therefore we have induction and restriction functors $\widetilde F$ and $\widetilde H$ going between $\Rep(\widetilde Q)$ and $\Rep(Q)$.

It now follows from Corollary \ref{Cor1} and Lemma \ref{Lem2} that there are natural isomorphisms $F\cong\widetilde H$ and $H\cong\widetilde F$.
\begin{Cor}
The ii-indecomposables for $(Q,\a)$ are given up to isomorphism by the images under $H$ of indecomposable $\widetilde Q$-representations, and $H(Y)\cong H(Z)$ if and only if $Y\cong {}^{\tilde\a^r}Z$ for some $r$.
\end{Cor}

\textit{Proof.} If $Y$ is an indecomposable for $\widetilde Q$, then $H(Y)$ is an ii-indecomposable for $(Q,\a)$. Conversely, suppose that $X$ is an indecomposable for $Q$ and let $Y$ be an indecomposable summand of $F(X)$. Then $X$ is a direct summand of $H(Y)$ and so every ii-indecomposable for $(Q,\a)$ is obtained up to isomorphism. \eop

We recall the construction of \cite{RR} for the semi-simple subalgebra of $kQ$. Let $\bi$ be a vertex orbit of $Q$ of size $d$ and let $i\in\bi$. Consider
$$A_\bi=k\varepsilon_i\times k\varepsilon_{\a(i)}\times\cdots\times k\varepsilon_{\a^{d-1}(i)}.$$
Using the  notation of \cite{RR}, we have an algebra isomorphism
$$\prod_{\mu=1}^{n/d}\mathbb M\,(d,k)_{(\bi,\mu)}\to A_\bi\#\langle\a\rangle, \qquad
E(\bi,\mu)_{pq}\mapsto\frac{d}{n}\sum_{j=1}^{n/d}\zeta^{dj\mu}\a^{-p}\varepsilon_i\a^{q+dj}$$
where $\zeta$ is a fixed primitive $n$-th root of unity and the $E(\bi,\mu)_{pq}$ are the elementary matrices. Let $\tilde\a$ be the generator for the dual group such that $\tilde\a(\a)=\zeta$. Then
$$\tilde\a\big(E(\bi,\mu)_{00}\big)=\tilde\a\Big(\frac{d}{n}\sum_j\zeta^{dj\mu}\varepsilon_i\a^{dj}\Big)=\frac{d}{n}\sum_j\zeta^{dj(\mu+1)}\varepsilon_i\a^{dj}=E(\bi,\mu+1)_{00}.$$
Write $E_\bi:=\sum_\mu E(\bi,\mu)_{00}$ and set $B_\bi=E_\bi A_\bi\#\langle\a\rangle E_\bi$. The algebras $A_\bi\#\langle\a\rangle$ and $B_\bi$ are then Morita equivalent, and since $E_\bi$ is fixed by $\tilde\a$ we have an induced action of the dual group on $B_\bi$.

In general, we repeat this for all vertex orbits $\bi$ and let $E=\sum_\bi E_\bi$ be the corresponding idempotent of $kQ\#\langle\a\rangle$, fixed by $\tilde\a$. Then $EkQ\#\langle\a\rangle E$ is isomorphic to the path algebra $k\widetilde Q$. The vertices $\widetilde\sI$ of $\widetilde Q$ are thus indexed by pairs $(\bi,\mu)$, where $\bi\in\bI$ and $\mu\in\bZ/(n/d_\bi)\bZ$.

We are now in a position to complete the proof of Theorem \ref{MainThm} --- namely that given any root $\alpha\in\Delta(\Gamma)_+$, there exists an indecomposable $\widetilde Q$-representation $Y$ such that $f(\underline\dim\, H(Y))=\alpha$. It is clearly enough to show that the map $h:\bZ\widetilde\sI\to\bZ\sI\to\bZ\bI$ induced by $H$ and $f$ maps $\Delta(\widetilde Q)_+$ onto $\Delta(\Gamma)_+$. Alos, the above description of the Morita equivalence implies that $h$ is given by
$$h(\beta)_\bi:=\sum_{\mu}\beta_{(\bi,\mu)}.$$ 

We first prove an analogous result to Lemma \ref{Lem}.
\begin{Lem}
Let $\tilde s_\bi:=\prod_\mu \tilde r_{(\bi,\mu)}\in W(\widetilde Q)$. We have the following relations.
\begin{enumerate}
\item $(h(\beta),e_\bi)_{\Gamma}=d_\bi\sum_\mu(\beta,e_{(\bi,\mu)})_{\widetilde Q}$;
\item $h(\tilde s_\bi(\beta))=r_\bi(h(\beta))$;
\item The map $r_\bi\mapsto \tilde s_\bi$ induces an isomorphism $W(\Gamma)\cong C_{\tilde\a}(W(\widetilde Q))$.
\end{enumerate}
\end{Lem}

\textit{Proof.} We note that if $\a^m$ fixes vertices $i$ and $j$, then we can find a basis of arrows $i\to j$ in $kQ$ with respect to which $\a^m$ acts diagonally (c.f. \cite{RR}). We now consider each of these basis elements separately. Let $\rho:i\to j$ be such a basis element and let $t$ be the lowest common multiple of $d_\bi$ and $d_\bj$. Then $\a^t(\rho)=\zeta^{rt}\rho$ for some $r$ and the arrows $$\rho,\a(\rho),\ldots,\a^{t-1}(\rho)$$ are linearly independent in $kQ$. From \cite{RR} we know that we get an arrow $(\bi,\mu)\to(\bj,\nu)$ in $k\widetilde Q$ from $\rho$ if and only if $\mu\equiv\nu+r\bmod n/t$, and so for fixed $\nu$ there are $t/d_\bi$ solutions for $\mu$. In particular, if $(\tilde a_{(\bi,\mu)(\bj,\nu)})$ is the GCM for $\widetilde Q$, then
$$d_\bi\sum_\mu \tilde a_{(\bi,\mu)(\bj,\nu)} = \sum_{v,w}a_{i_vj_w}=b_{\bi\bj} \qquad\textrm{for any }\nu$$
and hence
$$(h(\beta),e_\bi)_\Gamma = \sum_{\bi,\bj}b_{\bi\bj}h(\beta)_\bj = d_\bi\sum_\mu\sum_{(\bj,\nu)}\tilde a_{(\bi,\mu)(\bj,\nu)}\beta_{(\bj,\nu)} = d_\bi\sum_\mu(\beta,e_{(\bi,\mu)})_{\widetilde Q}.$$

The element $\tilde s_\bi$ is well defined since $\tilde\a$ is again admissible, and as the bilinear form $(-,-)_{\widetilde Q}$ is $\tilde\a$-invariant, $\tilde s_\bi$ commutes with the action of $\tilde\a$. Therefore
$$h(\tilde s_\bi(\beta))=h(\beta)-\sum_\mu(\beta,e_{(\bi,\mu)})_{\widetilde Q}e_\bi
=h(\beta)-\frac{1}{d_\bi}(h(\beta),e_\bi)_\Gamma e_\bi=r_\bi(h(\beta)).$$

As in the proof of Lemma \ref{Lem}, induction on length shows that $C_{\tilde\a}(W(\widetilde Q))$ is generated by the $\tilde s_\bi$. Finally, suppose that $r_{\bi_1}\cdots r_{\bi_m}=1\in W(\Gamma)$ and consider $\tilde s_{\bi_1}\cdots\tilde s_{\bi_m}\in C_{\tilde\a}(W(\widetilde Q))$. Then for any vertex $(\bj,\nu)$, $\tilde s_{\bi_1}\cdots\tilde s_{\bi_m}(e_{(\bj,\nu)})>0$. By the characterisation of the length of an element in the Weyl group, we must have that $\ell(\tilde s_{\bi_1}\cdots\tilde s_{\bi_m})=1$. It follows immediately that the $r_\bi$ and $\tilde s_\bi$ must satisfy the same relations. \eop

\begin{Prop}\label{Prop2}
The map $\beta\mapsto h(\beta)$ sends $\Delta(\widetilde Q)_+$ onto $\Delta(\Gamma)_+$. Moreover, if $\alpha\in\Delta(\Gamma)_+$ is real, then there is a unique $\tilde\a$-orbit of roots mapping to $\alpha$, all of which are real.
\end{Prop}

\textit{Proof.} Since the dimension vector of any ii-indecomposable must be a positive root for $\Gamma$, we know that $\beta\mapsto h(\beta)$ sends $\Delta(\widetilde Q)_+$ into $\Delta(\Gamma)_+$. (We could also prove this directly.) To show surjectivity, we first construct  preimages for the fundamental roots by adapting the proof of Lemma 5.3 in \cite{Kac4}.

Let $\alpha\in F_\Gamma$ and consider the set $\{\beta\in\Delta(\widetilde Q)_+\mid h(\beta)\leq\alpha\}$. Since this set is finite and non-empty, we can take an element $\beta$ of maximal height. Suppose $h(\beta)_\bi<\alpha_\bi$. Then for any $\mu$, $h(\beta+e_{(\bi,\mu)})=h(\beta)+e_\bi\leq\alpha$. By the maximality of $\beta$, $\beta+e_{(\bi,\mu)}$ cannot be a root and so $(\beta,e_{(\bi,\mu)})_{\widetilde Q}\geq 0$ (\cite{Kac4}, Corollary 3.6). Thus $(h(\beta),e_\bi)_\Gamma\geq0$ as well. In particular, $h(\beta)$ and $\alpha$ must have the same support, for otherwise we can find such a vertex $(\bi,\mu)$ adjacent to the support of $\beta$, and so $(\beta,e_{(\bi,\mu)})<0$. Contradiction.

We may assume that $\mathrm{supp}(\alpha)=\Gamma$. Let $S:=\{\bi\mid h(\beta)_\bi=\alpha_\bi\}$. If $S$ is the empty set, then $\beta+e_{(\bi,\mu)}$ is not a root for any vertex of $\widetilde Q$, and so the connected component of $\widetilde Q$ in which $\beta$ lies is Dynkin (\cite{Kac4}, Proposition 4.9). Therefore $\widetilde Q$ must be a disjoint union of copies of this Dynkin quiver, all in a single $\tilde\a$ orbit, so in particular $k\widetilde Q$ is representation finite. This implies that $kQ$ is representation finite (\cite{RR}, Theorem 1.3) and that $\Gamma$ is a connected Dynkin diagram. Contradiction, since $\alpha$ was assumed tobe imaginary.

Thus $S$ is non-empty so take any connected component $T$ of $\Gamma-S$ and write $\gamma$ for the restriction of $h(\beta)$ to $T$. Then for all vertices $\bj\in T$,
$$(\gamma,e_\bj)_T\geq(h(\beta),e_\bj)_\Gamma\geq 0.$$
Moreover, there exists a vertex $\bj\in T$ adjacent to $S$, and so $(\gamma,e_\bj)_T>0$. Therefore $T$ is Dynkin (\cite{Kac4}, Corollary 4.3).

Conversely, let $\gamma'$ be the restriction of $\alpha-h(\beta)$ to $T$. Note that $\gamma'$ has support the whole of $T$. Then, for any vertex $\bj\in T$,
$$(\gamma',e_\bj)_T=(\alpha-h(\beta),e_\bj)_\Gamma=(\alpha,e_\bj)_\Gamma-(h(\beta),e_\bj)_\Gamma\leq 0.$$
Hence $T$ is not Dynkin (\cite{Kac4}, Theorem 4.3). Contradiction. Therefore $S=\Gamma$ and $h(\beta)=\alpha$.

Clearly every simple root $e_\bi$ of $\Gamma$ lies in the image of $h$, and so the correspondence between the Weyl groups proves that $h$ is surjective.

Finally, let $\alpha$ be a real root for $\Gamma$ and let $\beta$ be a root for $\widetilde Q$ such that $h(\beta)=\alpha$. Then we can find an element $w'\in W(\Gamma)$ such that $w'(\alpha)$ is simple, say equal to $e_\bi$. If $w$ is the corresponding element in $C_{\tilde\a}(W(\widetilde Q))$, then $w(\beta)$ must also be simple, equal to some $e_{(\bi,\mu)}$. Therefore $\beta$ is real and uniquely determined up to an $\tilde\a$-orbit. \eop

Theorem \ref{MainThm} now follows immediately from Proposition \ref{Prop2}.

Using this proposition, we can exhibit a counter-example to the converse of Part 2 of Theorem \ref{MainThm} --- that is, an imaginary dimension vector  $\alpha\in\Delta(\Gamma)_+$ such that the corresponding ii-indecomposable is unique up to isomorphism. Namely, we consider
$$(Q,\a)\qquad\vcenter
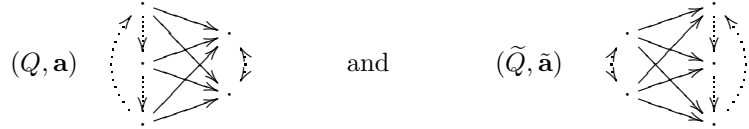
{\xymatrix@R=1pt{\cdot\ar[dr]\ar[dddr]\ar@{.>}[dd]\ar@/_1pc/@{<.}[dddd]\\&\cdot\ar@/^/@{<.>}[dd]\\\cdot\ar[ur]\ar[dr]
\ar@{.>}[dd]\\&\cdot\\\cdot\ar[uuur]\ar[ur]}} \qquad\qquad\textrm{and}\qquad\qquad(\widetilde Q,\tilde\a)\qquad
\vcenter{\xymatrix@R=1pt{&\cdot\ar@{.>}[dd]\ar@/^1pc/@{<.}[dddd]\\\cdot\ar[ur]\ar[dr]\ar[dddr]\ar@/_/@{<.>}[dd]\\&\cdot\ar@{.>}[dd]\\\cdot\ar[uuur]\ar[ur]\ar[dr]\\&\cdot}}$$
so that $\Gamma$ is the valued graph $\xymatrix@1{\cdot\ar@{-}[r]^{(3,2)}&\cdot}$. Then for $\alpha=(1,1)$ a fundamental root for $\Gamma$, there is a unique $\tilde\a$-orbit of roots for $\widetilde Q$ mapping to $\alpha$, all of which are real, and hence there is a unique ii-indecomposable of dimension vector $f^{-1}(\alpha)$.

\section{Tame Quivers and their Automorphisms}

We are interested in studying the tame quivers and their possible automorphisms. To this end, we first note that $\a$ commutes with the Auslander-Reiten translate $\tau$ (Lemma 4.1 in \cite{RR}). In particular, $\a$ acts on the preprojective and the preinjective components of the Auslander-Reiten quiver.

Now consider the special case when $Q$ is a tame quiver. Then the positive imaginary roots of $Q$ are generated by a single element $\delta$, and the defect of an indecomposable representation $X$ is defined to be
$$\mathrm{defect}(X):=\langle\delta,\underline\dim\, X\rangle_Q,$$
where $\langle-,-\rangle_Q$ is the Euler form for $Q$ given by
$$\langle\alpha,\beta\rangle_Q:=\sum_{i\in\sI}\alpha_i\beta_i-\sum_{\rho:i\to j}\alpha_i\beta_j.$$

The indecomposable regular representations are precisely those with defect 0, and $\tau$ acts with finite period on each connected regular component (tube) of the Auslander-Reiten quiver. Each regular indecomposable representation is regular serial, and is uniquely determined up to isomorphism by its regular top and regular length. Therefore the actions of both $\a$ and $\tau$ are determined by their actions on the regular simples.

We illustrate this with a couple of examples. Let $Q$ be the quiver $\widetilde D_4$
$$
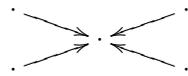
\xymatrix@R=1pt{\cdot\ar[dr]&&\cdot\ar[dl]\\&\cdot\\\cdot\ar[ur]&&\cdot\ar[ul]}$$
Then there are three tubes of period 2 with regular simples
\begin{alignat*}{3}
E_0&=\vcenter{\xymatrix@R=0pt@C=5pt{*+[d]{1}\ar@{-}[dr]&&*+[d]{0}\ar@{-}[dl]\\&1\\*+[u]{1}\ar@{-}[ur]&&*+[u]{0}\ar@{-}[ul]}}&\qquad
 E'_0&=\vcenter{\xymatrix@R=1pt@C=5pt{*+[d]{1}\ar@{-}[dr]&&*+[d]{1}\ar@{-}[dl]\\&1\\*+[u]{0}\ar@{-}[ur]&&*+[u]{0}\ar@{-}[ul]}}&\qquad
 E''_0&=\vcenter{\xymatrix@R=1pt@C=5pt{*+[d]{1}\ar@{-}[dr]&&*+[d]{0}\ar@{-}[dl]\\&1\\*+[u]{0}\ar@{-}[ur]&&*+[u]{1}\ar@{-}[ul]}}\\
E_1&=\vcenter{\xymatrix@R=1pt@C=5pt{*+[d]{0}\ar@{-}[dr]&&*+[d]{1}\ar@{-}[dl]\\&0\\*+[u]{1}\ar@{-}[ur]&&*+[u]{1}\ar@{-}[ul]}}&\qquad
 E'_1&=\vcenter{\xymatrix@R=1pt@C=5pt{*+[d]{0}\ar@{-}[dr]&&*+[d]{0}\ar@{-}[dl]\\&1\\*+[u]{1}\ar@{-}[ur]&&*+[u]{1}\ar@{-}[ul]}}&\qquad
 E''_1&=\vcenter{\xymatrix@R=1pt@C=5pt{*+[d]{0}\ar@{-}[dr]&&*+[d]{1}\ar@{-}[dl]\\&1\\*+[u]{1}\ar@{-}[ur]&&*+[u]{0}\ar@{-}[ul]}}
\end{alignat*}
as well as the tubes of period 1, indexed by $k-\{0,1\}$, with regular simples
$$T(\lambda)=\vcenter{\xymatrix@R=2pt{1\ar@{-}[dr]^{\tbinom{1}{0}}&&1\ar@{-}[dl]_{\tbinom{1}{1}}\\&2\\1\ar@{-}[ur]_{\tbinom{0}{1}}&&1\ar@{-}[ul]^
{\tbinom{1}{\lambda}}}}$$

If we let $\a$ act as the automorphism
$$
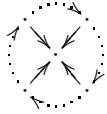
\xymatrix@=7pt@M=1pt{\cdot\ar@/^/@{.>}[rr]\ar[dr]&&\cdot\ar@/^/@{.>}[dd]\ar[dl]\\&\cdot\\\cdot\ar@/^/@{.>}[uu]\ar[ur]&&\cdot\ar@/^/@{.>}[ll]\ar[ul]}$$
so that $\Gamma$ is the valued graph $\widetilde A_{11}$ $\xymatrix@1{\cdot\ar@{-}[r]^{(4,1)}&\cdot}$, then $\a$ acts on the regular simples (up to isomorphism) as
$$E_0\mapsto E'_0\mapsto E_1\mapsto E'_1, \qquad E''_0\mapsto E''_1, \qquad T(\lambda)\mapsto T(\tfrac{\lambda}{\lambda-1}).$$
Therefore we have ii-indecomposables
$$E_0\oplus E'_0\oplus E_1\oplus E'_1, \quad E''_0\oplus E''_1, \quad T(2), \quad T(\lambda)\oplus T(\tfrac{\lambda}{\lambda-1}) \textrm{ for } \lambda\neq 0,1,2.$$

Alternatively, if $\a$ acts as the automorphism
$$
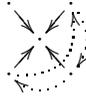
\xymatrix@=7pt@M=1pt{\cdot\ar[dr]&&\cdot\ar@/^/@{.>}[dd]\ar[dl]\\&\cdot\\\cdot\ar@/_1pc/@{.>}[uurr]\ar[ur]&&\cdot\ar@/^/@{.>}[ll]\ar[ul]}$$
then $\Gamma$ is the valued graph $\widetilde G_{21}$ $\xymatrix@1@C=15pt{\cdot\ar@{-}[r]&\cdot\ar@{-}[r]^{(1,3)}&\cdot}$ and we have ii-indecomposables
\begin{gather*}
E_0\oplus E'_0\oplus E''_0, \quad E_1\oplus E'_1\oplus E''_1, \quad T(-\omega), \quad T(-\omega^2)\\
T(\lambda)\oplus T(\tfrac{1}{1-\lambda})\oplus T(\tfrac{\lambda-1}{\lambda})\textrm{ for }\lambda\neq 0,1,\omega,\omega^2,
\end{gather*}
where $\omega$ is a primitive cube root of unity.

\section{Representations of Species over Finite Fields}

Let $\Gamma$ be a valued quiver with symmetrisable GCM $C=D^{-1}B$. We fix a finite field $k=\bF_q$, an algebraic closure $\bar k$ and denote by $\bF_{q^r}$ the unique extension of $k$ of degree $r$ inside $\bar k$. The $k$-species $S$ of $\Gamma$ is given by the field $\bF_{q^{d_i}}$ at vertex $i$ and the $\bF_{q^{d_j}}$-$\bF_{q^{d_i}}$-bimodule $\bF_{q^{b_{ij}}}$ for each arrow $i\to j$ in $\Gamma$.

A representation $X$ of $S$ is given by an $\bF_{q^{d_i}}$-vector space $X_i$ for each vertex $i$ and an $\bF_{q^{d_j}}$-linear map
$$\theta_{ij}:\bF_{q^{b_{ij}}}\otimes_{\bF_{q^{d_i}}}X_i\to X_j$$
for each arrow $i\to j$. Therefore, the category of $S$-representations is equivalent to the category of finite demensional modules for the tensor algebra $\Lambda:=T(A_0,A_1)$, where $A_0=\prod_i\bF_{q^{d_i}}$ is a semisimple algebra and $A_1=\coprod_{i\to j}\bF_{q^{b_{ij}}}$ is an $A_0$-bimodule.

In \cite{DR}, Dlab and Ringel studied the representations of a $k$-species associated to a valued Dynkin or tame quiver $\Gamma$. For $k$ a finite field, the number of isomorphism classes of representations was studied for general $\Gamma$ (without oriented cycles) by Hua \cite{Hua2} and a partial analogue of Kac's Theorem was proved. This was later completed by Deng and Xiao \cite{DX} using the Ringel-Hall algebra.

\begin{Thm}[Hua, Deng-Xiao]
Let $\Gamma$ be a valued quiver, $k$ a finite field and $S$ the $k$-species of $\Gamma$. Then
\begin{enumerate}
\item the dimension vectors of the indecomposable $S$-representatons are precisely the positive roots of the symmetrisable Kac-Moody Lie algebra $\mathfrak g(\Gamma)$;
\item if $\alpha$ is real, there is a unique indecomposable of dimension vector $\alpha$ up to isomorphism.
\end{enumerate}
\end{Thm}

Suppose now that $K/k$ is a finite field extension. Then the Galois group $G:=\mathrm{Gal}(K,k)$ acts as $k$-algebra automorphisms on $K\otimes_k\Lambda$, say generated by $g:a\otimes\lambda\mapsto a^q\otimes\lambda$. Thus for any $K\otimes\Lambda$-module $Y$ we can form a new module ${}^gY$ by taking the same underlying $k$-vector space with the new action $p\cdot y:=g^{-1}(p)y$.
\begin{Prop}\label{Prop3}
There is a bijection between the isomorphism classes of $\Lambda$-modules and the isomorphism classes of $K\otimes_k\Lambda$-modules $Y$ satisfying ${}^gY\cong Y$.
\end{Prop}

The proof is based upon the following series of lemmas (c.f. \cite{KR}).
\begin{Lem}
Let $X$ be a $\Lambda$-module and $K/k$ a field extension of degree $r$. Then $K\otimes_kX|_k\cong X^r$ over $k$. In particular, two $\Lambda$-modules $X$ and $Y$ are isomorphic if and only if $K\otimes_kX\cong K\otimes_kY$.
\end{Lem}

\begin{Lem}\label{Lem3}
Let $X$ be an indecomposable $\Lambda$-module, $K:=\mathrm{End}(X)/\mathrm{rad}\,\mathrm{End}(X)$ and $r=[K:k]$. Then $K\otimes_kX$ is a direct sum of $r$ pairwise non-isomorphic indecomposables $Y_i$ with $\mathrm{End}(Y_i)/\mathrm{rad}\,\mathrm{End}(Y_i)=K$.
\end{Lem}

\textit{Proof.} Since $k$ is perfect, $\mathrm{End}(K\otimes_kX)/\mathrm{rad}\,\mathrm{End}(K\otimes_kX)$ is isomorphic to $K\otimes_k\mathrm{End}(X)/\mathrm{rad}\,\mathrm{End}(X)\cong K\otimes_kK\cong K^r$. \eop

\begin{Lem}
Let $K/k$ be a finite field extension with Galois group $G$. Let $Y$ be a $K\otimes_k\Lambda$-module such that ${}^gY\cong Y$ for all $g\in G$ and suppose that $G$ acts transitively on the indecomposable summands of $Y$ (up to isomorphism). Then any indecomposable summand $X$ of the $\Lambda$-module $Y|_k$ satisfies $K\otimes_kX\cong Y$.
\end{Lem}

\textit{Proof.}
We have $K\otimes_kY\cong\oplus_{g\in G}{}^gY\cong Y^r$, using that $K\otimes_kK\cong\oplus_{g\in G}{}^gK$ as $K$-bimodules. Now let $X$ be an indecomposable summand of $Y|_k$. Since ${}^g(K\otimes_kX)\cong K\otimes_kX$, the assumption on $Y$ implies that $K\otimes_kX\cong Y^s$ for some $s$. If $L$ is common extension of $K$ and $\mathrm{End}(X)/\mathrm{rad}\,\mathrm{End}(X)$, then on the one hand, $M\otimes_kX$ is isomorphic to the direct sum of pairwise non-isomorphic indecomposables, whereas on the other it is isomorphic to $M\otimes_KY^s$. Thus $s=1$. \eop

In fact, we now see that if $X$ is a $\Lambda$-indecomposable as in Lemma \ref{Lem3}, then $G$ must act transitively on the $Y_i$. Proposition \ref{Prop3} now follows immediately.

In particular, we can take $K/k$ to be a splitting field for each $\bF_{q^{d_i}}$ --- for example, $K=\bF_{q^t}$ where $t$ is the lowest common multiple of the $d_i$. Then the algebra $K\otimes_k\Lambda$ is isomorphic to the path algebra of a quiver $Q$. We show that the action of the Galois group can now be thought of as an admissible quiver automorphism combined with the $G$-action on $KQ$ coming from the identification $KQ\cong K\otimes_kkQ$.

We describe this construction in the case when $\Gamma$ has only two vertices; the general case following immediately. Let $\Gamma$ be the valued quiver
$$\xymatrix@1@C=40pt{i\ar[r]^{(c_{ji},c_{ij})} &j}$$
with $d_ic_{ij}=b_{ij}=d_jc_{ji}$. Let $t$ be the lowest common multiple of $d_i$ and $d_j$ and set $b_{ij}=bt$. The tensor algebra is therefore given by
$$\Lambda=\begin{pmatrix}\bF_{q^{d_i}} & 0\\\bF_{q^{bt}} & \bF_{q^{d_j}} \end{pmatrix}
\cong\begin{pmatrix}\bF_{q^{d_i}} & 0\\K^b & \bF_{q^{d_j}}\end{pmatrix}.$$
Now $K\otimes_k\Lambda\cong KQ$ for some quiver $Q$, and we can assume that $b=1$ since for general $b$ we just take $b$ copies of each arrow in $Q$.

By the Normal Basis Theorem \cite{Coh}, there exists $x\in\bF_{q^{d_i}}$ such that the elements $x^{q^r}$ for $0\leq r<d_i$ form a $k$-basis for $\bF_{q^{d_i}}$. We also have the isomorphism
\begin{equation}\label{eqtn}
\bF_{q^{d_i}}\otimes_k\bF_{q^{d_i}}\overset{\sim}{\to}\prod_{\mu=0}^{d_i-1}\bF_{q^{d_i}}, \qquad a\otimes b\mapsto (ab^{q^\mu})_\mu.\tag{*}
\end{equation}
Let $\varepsilon_{(i,\mu)}$ denote the unit for the $\mu$-th copy of $\bF_{q^{d_i}}$ and fix $a_r\in\bF_{q^{d_i}}$ such that $\sum_ra_r\otimes x^{q^r}$ corresponds to $\varepsilon_{(i,0)}$. Then $\sum_ra_r\otimes x^{q^{r-\mu}}$ corresponds to $\varepsilon_{(i,\mu)}$ and so
$$\sum_{r,s}a_ra_s\otimes x^{q^{r-\lambda}}x^{q^{s-\mu}}=\delta_{\lambda,\mu}\sum_ra_r\otimes x^{q^{r-\lambda}}.$$

Similarly we can find $y,b_s\in K$ such that the $y^{q^{sd_i}}$ give an $\bF_{q^{d_i}}$-basis for $K$ and $\sum_sb_sy^{q^{(s-\mu)d_i}}=\delta_{\mu,0}$. Then, for the isomorphism
$$K\otimes_kK\xrightarrow{\sim}\prod_{\lambda=0}^{t-1}K, \qquad a\otimes b\mapsto (ab^{q^\lambda})_\lambda,$$
the $\lambda$-th unit, which we denote by $\rho_{(\lambda)}$, corresponds to $\sum_{r,s}a_rb_s\otimes x^{q^{r-\lambda}}y^{q^{sd_i-\lambda}}$. It follows that
$$\rho_{(\lambda)}\varepsilon_{(i,\mu)}=\begin{cases}\rho_{(\lambda)} &\textrm{if }\lambda\equiv\mu\bmod d_i;\\0&\textrm{otherwise}. \end{cases}$$

We have analogous identities for the isomorphism $K\otimes_k\bF_{q^{d_j}}\xrightarrow{\sim}\prod_{\nu=0}^{d_j-1}K$ and the product $\varepsilon_{(j,\nu)}\rho_{(\lambda)}$. In conclusion, we have

\begin{Lem}
Let $\Gamma$ be a valued quiver with symmetrisable GCM $\ C=D^{-1}B$, $k$ a finite field and $\Lambda$ the corresponding $k$-algebra. Let $t$ be the lowest common multiple of the $d_i$. Then for $K/k$ an extension of degree $t$, the algebra $K\otimes\Lambda$ is isomorphic to the path algebra of a quiver $KQ$. The vertices of $Q$ are labelled by pairs $(i,\mu)$ for $i$ a vertex of $\Gamma$ and $0\leq\mu<d_i$ and there are $a_{ij}$ arrows $(i,\mu)\to(j,\nu)$ if and only if $i\to j$ in $\Gamma$ and $\mu\equiv\nu\bmod\mathrm{hcf}(d_i,d_j)$, where $a_{ij}:=b_{ij}/\mathrm{lcm}(d_i,d_j)$.
\end{Lem}

Let $G$ denote the Galois group $\textrm{Gal}(K,k)$. This acts naturally on $\Lambda$ as $k$-algebra automorphisms, and so we can extend this to get $K$-algebra automorphisms of $K\otimes\Lambda$. Let $\a$ be a generator for $G$, acting as $\a(a\otimes b):=a\otimes b^q$. In terms of the path algebra $KQ$, we get
$$\a(\varepsilon_{(i,\mu)})=\a\big(\sum_ra_r\otimes x^{q^{r-\mu}}\big)=\sum_ra_r\otimes x^{q^{r+1-\mu}}=\varepsilon_{(i,\mu-1)}.$$
Similarly, $\a(\varepsilon_{(j,\nu)})=\varepsilon_{(j,\nu-1)}$ and $\a(\rho_{(\lambda)})=\rho_{(\lambda-1)}$, and hence $\a$ acts on $KQ$ as an admissible quiver automorphism. We note that the valued graph associated to the pair $(Q,\a)$ is precisely the underlying graph of the valued quiver $\Gamma$.

On the other hand, we have the standard $G$-action on $KQ$ coming from the isomorphism $KQ\cong K\otimes_kkQ$. Let $\tau$ be a generator for this action, where for example $\tau(a\varepsilon_{(i,\mu)})=a^q\varepsilon_{(i,\mu)}$. Then $\tau$ acts on $K\otimes\Lambda$ as $a\otimes b\mapsto a^q\otimes b^q$. In particular, the $k$-algebra automorphism $\a^{-1}\tau$ of $K\otimes\Lambda$ sends $a\otimes b$ to $a^q\otimes b$ and so we recover our original $G$ action. Thus we deduce the following proposition.

\begin{Prop}
There is a bijection between the isomorphism classes of $\Lambda$-modules of dimension vector $\alpha$ and representations $Y$ of $KQ$ of dimension vector $f^{-1}(\alpha)$ such that ${}^\a Y\cong{}^\tau Y$, where $f$ is again the canonical map from the $\a$-fixed points of the root lattice for $Q$ to the root lattice for $\Gamma$.
\end{Prop}

We can now offer a more representation-theoretic poof of the generalisation of Kac's Theorem to representations of species over finite fields (c.f. \cite{Hua2, DX}). We denote by $I_\Gamma(\alpha,q)$ the number of isomorphism classes of indecomposable representations of dimension vector $\alpha$ for the $\bF_q$-species of $\Gamma$.

\begin{Thm}[Hua \cite{Hua2}]
The numbers $I_\Gamma(\alpha,q)$ are polynomial in $q$ with rational coefficients, independent of the orientation of $\Gamma$.
\end{Thm}

\begin{Thm}
The polynomials $I_\Gamma(\alpha,q)$ are non-zero if and only if $\alpha$ is a positive root of $\mathfrak g(\Gamma)$. Moreover, if $\alpha$ is a real root of $\mathfrak g(\Gamma)$ then $I_\Gamma(\alpha,q)=1$.
\end{Thm}

\textit{Proof.} Let $\Gamma$ have symmetrisable GCM $C=D^{-1}B$, let $S$ be the $k$-species of $\Gamma$ and write $\Lambda$ for the corresponding tensor algebra. Set $t$ to be the lowest common multiple of the $d_i$ and let $K/k$ be a field extension of degree $t$. Then $K\otimes\Lambda$ is isomorphic to a path algebra $KQ$ and as before, we can consider the two different actions of $G=\textrm{Gal}(K,k)$ on $KQ$ generated by $\a$ and $\tau$.

We know that the isomorphism classes of indecomposable $\Lambda$-modules are in bijection with the isomorphism classes of representations $Y$ for $KQ$ such that ${}^\a Y\cong{}^\tau Y$ and $Y$ is not the proper direct sum of two such representations. In particular, every ii-indecomposable $Y$ for $(Q,\a)$ of dimension vector $f^{-1}(\alpha)$ defined over $k$ gives rise to a $\Lambda$-indecomposable of dimension vector $\alpha$.

Let $X$ be a $\Lambda$-indecomposable of dimension vector $\alpha$ and let $Y$ be an indecomposable summand of $K\otimes X$. Setting $\underline\dim\,Y=\beta$, we see that $\underline\dim\,K\otimes X=f^{-1}(\alpha)=r\sigma(\beta)$ for some $r$. Thus if $\beta$ is imaginary root of $Q$, then $\alpha=rf(\sigma(\beta))$ is an imaginary root of $\Gamma$.

Conversely, if $\beta$ is real, then $Y$ is defined over $k$ and so fixed by $\tau$. Therefore $K\otimes X$ is an ii-indecomposable for $(Q,\a)$. Hence $r=1$ and $\alpha$ is a root of $\Gamma$. This shows that all indecomposable $\Lambda$-modules have dimension vector a root of $\Gamma$ and hence the polynomial $I_\Gamma(\alpha,q)$ is non-zero unless $\alpha\in\Delta(\mathfrak g)_+$ (C.f. the argument in \cite{Hua2}).

Moreover, if $\alpha$ is real, then the highest common factor of the $\alpha_i$ is 1 and so $r=1$. Therefore $\beta$ is real and unique up to $\a$-orbit. Since there is a unique isomorphism class of ii-indecomposables for $(Q,\a)$ of dimension vector $f^{-1}(\alpha)$, we deduce that $I_\Gamma(\alpha,q)=1$. 

Finally, suppose that $\alpha\in\Delta(\Gamma)_+$. Then there exists an ii-indecomposable $X$ for $(Q,\a)$ of dimension vector $f^{-1}(\alpha)$ over $\overline\bF_p$, where $p$ is coprime to $t$. This $X$ must be defined over some finite field $k=\bF_q$ and so corresponds to an indecomposable for the $k$-species of $\Gamma$ of dimension vector $\alpha$. Therefore the polynomial $I_\Gamma(\alpha,q)$ cannot be zero. \eop

\end{document}